\documentclass{gtart}

\def\ifplaintex{\expandafter\ifx\csname documentclass\endcsname\relax}

\ifplaintex 
\else
\fi

\expandafter\ifx\csname beginpicture\endcsname\relax
\expandafter\ifx\csname documentclass\endcsname\relax
\input pictex \else
\input prepictex \input pictex \input postpictex \fi\fi

\def\gt{{\mathsurround=0pt\it $\cal G\mskip-2mu$eometry \&\ 
$\cal T\!\!$opology}}        

\def\gtp{{\mathsurround=0pt\it $\cal G\mskip-2mu$eometry \&\ 
$\cal T\!\!$opology $\cal P\!$ublications}}  

\def\lognumber#1{\def\thelognumber{#1}}
\def\volumenumber#1{\def\thevolumenumber{#1}}
\def\papernumber#1{\def\thepapernumber{#1}}
\def\volumeyear#1{\def\thevolumeyear{#1}}

\def\pagenumbers#1#2{\def\startpage{#1}\def\finishpage{#2}}
\def\published#1{\def\publishdate{#1}}
\def\proposed#1{\def\theproposer{#1}}
\def\seconded#1{\def\theseconders{#1}}
\def\received#1{\def\receiveddate{#1}}

\def\accepted#1{\def\accepteddate{#1}}

\long\def\asciiabstract#1{\long\def\theasciiabstract{#1}}
\def\asciikeywords#1{\def\theasciikeywords{#1}}

\def\shorttitle#1{\def\theshorttitle{#1}}

\let\\\par\let\thelognumber\relax
\let\thevolumenumber\relax\let\thepapernumber\relax
\let\thevolumeyear\relax\let\thesamplenumber\relax\let\startpage\relax
\let\finishpage\relax\let\publishdate\relax\let\receiveddate\relax
\let\reviseddate\relax\let\accepteddate\relax\let\theasciititle\relax
\let\theasciiauthors\relax
\let\theasciiabstract\relax\let\theasciikeywords\relax
\let\theasciiemail\relax\let\theshortauthors\relax\let\theshorttitle\relax

\long\def\maketitlep{   

\count0=\startpage

\gt\hfill      
\beginpicture
\setcoordinatesystem units <0.33truein, 0.33truein> point at 2.2 0.9
\setplotsymbol ({$\cal G$})
\plotsymbolspacing=9truept
\circulararc 315 degrees from 0 1 center at 0 0
\setplotsymbol ({$\cal T$})
\circulararc 315 degrees from 1 -1 center at 1 0
\endpicture
\break
{\small\ifx\thesamplenumber\relax 
Volume \else Sample
\fi\thevolumenumber\ (\thevolumeyear)
\startpage--\finishpage\nl
Published: \publishdate}
\vglue 0.5truein plus 0.4fil minus 0.1truein

{\parskip=0pt\leftskip 0pt plus 1fil\def\\{\par\smallskip}{\ifplaintex\large
\else\Large\fi\bf\thetitle}\par\medskip}   

\vglue 0pt plus 0.1fil 

{\parskip=0pt\leftskip 0pt plus 1fil\def\\{\par}{\sc\theauthors}
\par\medskip}

\vglue 0pt plus 0.1fil 

{\small\parskip=0pt\let\newline\\
{\leftskip 0pt plus 1fil\def\\{\par}{\sl\theaddress}\par}
\expandafter\ifx\theemail\relax    
\relax\else\vglue 5pt plus 0.02fil minus 2pt\def\\{\stdspace{\rm 
and}\stdspace} 
\cl{Email:\stdspace\tt\theemail}\fi
\ifx\theurl\relax                  
\relax\else\vglue 5pt plus 0.02fil minus 2pt\def\\{\stdspace{\rm 
and}\stdspace}
\cl{URL:\stdspace\tt\theurl}\fi\par}

\vglue 7pt plus 0.3fil minus 3pt

{\bf Abstract}
\vglue 5pt plus 0.1fil minus 2pt

\theabstract

\vglue 7pt plus 0.3fil minus 3pt

{\bf AMS Classification numbers}\quad Primary:\quad \theprimaryclass

Secondary:\quad \thesecondaryclass

\vglue 5pt plus 0.3fil minus 2pt

{\bf Keywords}\quad \thekeywords

\vglue 10pt plus 0.5fil minus 5pt

{\small  Proposed: \theproposer\hfill Received: \receiveddate\nl
Seconded: \theseconders\hfill 
\ifx\reviseddate\relax                         
Accepted: \accepteddate                        
\else
Revised: \reviseddate                          
\fi}
\eject
}       

\let\maketitlepage\maketitlep
\let\maketitle\maketitlepage

\font\phead=cmsl9 scaled 950
\font=cmsl9 scaled 1050
\font\pnum=cmbx10 scaled 913
\font\lnum=cmbx10 
\font\pfoot=cmsl9 scaled 950
\font=cmsl9 scaled 1050
\ifplaintex
\headline{\vbox to 0pt{\vskip -4.5mm\line{\small\phead\ifnum
\count0=\startpage ISSN 1364-0380 (on line)
1465-3060 (printed) \hfill {\pnum\folio}\else\ifodd\count0\def\\{ }%
\ifx\theshorttitle\relax\thetitle\else\theshorttitle\fi\hfill{\pnum\folio}
\else\def\\{ and }{\pnum\folio}\hfill\ifx\theshortauthors\relax\theauthors
\else\theshortauthors\fi\fi\fi}\vss}}
\footline{\vbox to 0pt{\vglue 0mm\line{\small\pfoot\ifnum\count0=\startpage
\copyright\ \gtp\hfill\else
\gt, Volume \thevolumenumber\ (\thevolumeyear)\hfill\fi}\vss
}}
\else
\makeatletter
\def\@oddhead{{\small\lhead\ifnum\count0=\startpage ISSN 1364-0380 (on line)
1465-3060 (printed) \hfill {\lnum\number\count0}\else\ifodd\count0
\def\\{ }\ifx\theshorttitle\relax \thetitle \else\theshorttitle\fi\hfill
{\lnum\number\count0}\else\def\\{ and }{\lnum\number\count0}
\hfill\ifx\theshortauthors\relax 
\theauthors\else\theshortauthors\fi\fi\fi}}\def\@evenhead{\@oddhead}
\def\@oddfoot{\small\lfoot\ifnum\count0=\startpage\copyright\ \gtp\hfill\else
\gt, Volume \thevolumenumber\ (\thevolumeyear)\hfill\fi}
\def\@evenfoot{\@oddfoot}
\makeatother
\fi

\newwrite\gtoutfile
\long\gdef\makeheadfile{  
{\def\\{, }\def\s{ }
\immediate\openout\gtoutfile head.xxx
\immediate\write\gtoutfile{To: math@arxiv.org}
\immediate\write\gtoutfile{Subject: put or rep NNNNN:pppp}
\immediate\write\gtoutfile{--text follows this line--}
\immediate\write\gtoutfile{Proxy-for: \ifx\theasciiauthors\relax
\theauthors\else\theasciiauthors\fi\s<\ifx\theasciiemail\relax\theemail\else\theasciiemail\fi>}
\immediate\write\gtoutfile{\noexpand\\}
\immediate\write\gtoutfile{Authors: \ifx\theasciiauthors\relax
\theauthors\else\theasciiauthors\fi}
{\def\\{ }\immediate\write\gtoutfile{Title: \ifx\theasciititle\relax
\thetitle\else\theasciititle\fi}}
\immediate\write\gtoutfile{Subj-class: GT or SG or MG etc}
\immediate\write\gtoutfile{MSC-class: \theprimaryclass\ifx\thesecondaryclass\relax\else, \thesecondaryclass\fi}
\immediate\write\gtoutfile{Journal-ref: Geom. Topol. \thevolumenumber
(\thevolumeyear) \startpage-\finishpage}
\immediate\write\gtoutfile{Comments: Published by Geometry and Topology at}
\immediate\write\gtoutfile{\s\s http://www.maths.warwick.ac.uk/gt/GTVol\thevolumenumber/paper\thepapernumber.abs.html}
\immediate\write\gtoutfile{\noexpand\\}
\immediate\write\gtoutfile{}
\ifx\theasciiabstract\relax
\immediate\write\gtoutfile{\theabstract}\else
\immediate\write\gtoutfile{\theasciiabstract}\fi
\immediate\write\gtoutfile{}
\immediate\write\gtoutfile{\noexpand\\}
\immediate\write\gtoutfile{}
\immediate\closeout\gtoutfile}}  

\def\maketitlepage{\maketitlep\makeheadfile}
\let\maketitle\maketitlepage

\lognumber{295}
\volumenumber{7}\papernumber{14}
\volumeyear{2003}
\pagenumbers{487}{510}
\received{9 December 2002}
\accepted{10 July 2003}
\published{29 July 2003}
\proposed{Gang Tian}
\seconded{John Morgan, Leonid Polterovich}

\usepackage{amsmath, amssymb}

\newtheorem{theorem}{Theorem}

\newtheorem{corollary}[theorem]{Corollary}

\newtheorem{lemma}[theorem]{Lemma}
\newtheorem{proposition}[theorem]{Proposition}

\theoremstyle{definition}
\newtheorem{definition}[theorem]{Definition}

\newtheorem{example}[theorem]{Example}

\begin{document}

\title{Precompactness of solutions to the Ricci flow\\in the absence of 
injectivity radius estimates}
\shorttitle{Precompactness of solutions to the Ricci flow}
\author{David Glickenstein}
\address{Department of Mathematics,
University of California, San Diego\\
9500 Gilman Drive,
La Jolla, CA 92093-0112, USA}

\email{glicken@math.ucsd.edu}

\asciiabstract{Consider a sequence of pointed n-dimensional complete Riemannian manifolds
{(M_i,g_i(t), O_i)} such that t in [0,T] are solutions to the Ricci flow and
g_i(t) have uniformly bounded curvatures and derivatives of curvatures. Richard
Hamilton showed that if the initial injectivity radii are uniformly bounded
below then there is a subsequence which converges to an n-dimensional solution
to the Ricci flow. We prove a generalization of this theorem where the initial
metrics may collapse. Without injectivity radius bounds we must allow for
convergence in the Gromov-Hausdorff sense to a space which is not a manifold
but only a metric space. We then look at the local geometry of the limit to
understand how it relates to the Ricci flow.}

\asciikeywords{Ricci flow, Gromov-Hausdorff convergence}

\begin{abstract}
Consider a sequence of pointed $n$--dimensional complete Riemannian
manifolds $\{(M_i,g_i(t), O_i)\}$ such that $t\in[0,T]$ are solutions
to the Ricci flow and $g_i(t)$ have uniformly bounded curvatures and
derivatives of curvatures. Richard Hamilton showed that if the initial
injectivity radii are uniformly bounded below then there is a
subsequence which converges to an $n$--dimensional solution to the Ricci
flow. We prove a generalization of this theorem where the initial
metrics may collapse. Without injectivity radius bounds we must allow
for convergence in the Gromov--Hausdorff sense to a space which is not
a manifold but only a metric space. We then look at the local geometry
of the limit to understand how it relates to the Ricci flow.
\end{abstract}




\maketitlepage

\section{Introduction}

Consider a solution to the Ricci flow $\left(  \mathcal{M},g\left(  t\right)
\right)  $ defined on a maximal time interval $t\in\left[  0,T_{M}\right)  .$
By Richard Hamilton's work \cite{hamilton:threemanifolds} the curvature must
go to infinity as $t\rightarrow T_{M}$, but we can look at the limit of
dilations of the singularity by taking a sequence of points $x_{i}%
\in\mathcal{M}$ and times $t_{i}\nearrow T_{M}$ such that
\[
\left|  Rm\left(  x_{i},t_{i}\right)  \right|  _{g_{i}\left(  t_{i}\right)
}=\sup_{x\in\mathcal{M}}\left|  Rm\left(  x,t_{i}\right)  \right|
_{g_{i}\left(  t_{i}\right)  }.
\]
We can then consider rescaled solutions to the Ricci flow such as%
\[
g_{i}\left(  t\right)  =\left|  Rm\left(  x_{i},t_{i}\right)  \right|
\,g\left(  t_{i}+\frac{t}{\left|  Rm\left(  x_{i},t_{i}\right)  \right|
}\right)
\]
so that we get a sequence of solutions to the Ricci flow on $\mathcal{M}$ with
curvature bounded above and below. We would like to look at the limit of such
metrics on $\mathcal{M}$ in order to study what happens as the Ricci flow
approaches a singularity. This is referred to as ``blowing up about a
singularity.'' Note two things about the blow up:

\begin{enumerate}
\item  On a fixed compact time interval the curvatures are bounded.

\item  The diameters may go to infinity.
\end{enumerate}

With these sorts of blow ups in mind Richard Hamilton proved a compactness
theorem for solutions to the Ricci flow in \cite{hamilton:compactness}. He
defines a complete marked solution to the Ricci flow $\left(  \mathcal{M}%
^{n},g\left(  t\right)  ,O,\mathcal{F}\right)  $ where $\mathcal{M}$ is an
$n$--dimensional $C^{\infty}$ manifold, $g\left(  t\right)  $ is a complete
Riemannian metric evolving by the Ricci flow on some closed time interval
$\left[  0,T\right]  ,$ $O\in\mathcal{M}$ is a basepoint and $\mathcal{F}$ is
a frame at $O$ orthonormal with respect to the metric $g\left(  0\right)  .$
The theorem then states that every sequence of complete marked solutions to
the Ricci flow $\left\{  \left(  \mathcal{M}_{i}^{n},g_{i}\left(  t\right)
,O_{i},\mathcal{F}_{i}\right)  \right\}  _{i\in\mathbb{N}}$ with uniformly
bounded curvature tensors and with injectivity radii of $g_{i}\left(
0\right)  $ at $O_{i}$ uniformly bounded below has a subsequence which
converges to a complete marked solution to the Ricci flow $\left(
\mathcal{M}_{\infty}^{n},g_{\infty}\left(  t\right)  ,O_{\infty}%
,\mathcal{F}_{\infty}\right)  .$ Convergence means that for every compact set
$K\subset\mathcal{M}_{\infty}$ and for $i$ large enough there are smooth
embeddings $\phi_{i}:K\rightarrow\mathcal{M}_{i}$ such that $\phi_{i}^{\ast
}\left(  g_{i}\right)  $ converges uniformly to $\left.  g_{\infty}\right|
_{K}.$

We choose the blow up about a singularity so that the curvature bound is
immediately satisfied, but the injectivity radius bound is not immediate. For
this reason much work has been done to prove that an injectivity radius bound
exists for solutions to the Ricci flow. The type I\ case was studied by
Hamilton in \cite{hamilton:singularities}. Also the case of type II with bumps
of curvature has been studied by Hamilton in \cite[Section 23]%
{hamilton:singularities} and Chow--Knopf--Lu in \cite{chowknopflu}. More
recently, Perelman \cite{perelman} proved an injectivity radius estimate for
finite time singularities on compact manifolds. In \cite[Section
2]{hamilton:compactness} Richard Hamilton suggests considering the case where
there is no injectivity bounds by looking at the work of Fukaya. This is our
main motivation. We also note that Carfora--Marzuoli \cite{carforamarzuoli}
first used Gromov--Hausdorff convergence as a tool to look at Ricci flow.

Much work has been done studying convergence of Riemannian manifolds, for
instance the work by Greene and Wu \cite{greenewu}, Peters \cite{peters}, and
Petersen \cite{petersen} for $C^{k}$ convergence. Hamilton notes, however,
that when we are concerned with the Ricci flow much of the hard work is wasted
and it is just as easy to work with basic principles. We follow this
philosophy. Gromov \cite{gromovenglish} proved that $n$--dimensional Riemannian
manifolds with Ricci curvature bounded below are precompact in the pointed
Gromov--Hausdorff topology. Much work has been done to try to understand the
collapsing case. Some noteworthy work is that of Cheeger--Gromov
\cite{cheegergromov1} and \cite{cheegergromov2} looking at F--structures,
Fukaya's fibration theorem for bounded curvature collapse in
\cite{fukayacollapse} and \cite{fukayacollapse2} generalized to lower
curvature bound by Yamaguchi \cite{yamaguchi}, Cheeger--Fukaya--Gromov's
development of nilpotent Killing structures in \cite{cheegerfukayagromov},
Shioya--Yamaguchi's study of three manifold collapse with lower curvature bound
\cite{shioyayamaguchi}, Cheeger--Rong's study of collapse with bounded covering
geometry in \cite{cheegerrong}, and work on collapsing with curvature pinched
positive by Petrunin--Rong--Tuschmann \cite{petruninrongtuschmann} and Rong
\cite{rong}.

Our work is modelled on Fukaya's study of the local structure of the limit of
Riemannian manifolds with a two sided bound on curvature in
\cite{fukayaboundary}. These results will be reviewed in Section
\ref{backgroundsection}. The following definition is inspired by Fukaya's
study of the local structure of limits. Actually, there is a more precise
definition that can be given inspired by Peter Petersen's definition of
convergence for the noncollapsing case in \cite{petersen}. In this case we can
define a structure on certain metric spaces that are controlled by curvature
bounds and how large our submersions are. In this paper we restrict to the
simpler case.

\begin{definition}
A sequence of marked $n$--dimensional Riemannian manifolds $\{(\mathcal{M}%
_{i}^{n},\allowbreak g_{i},\allowbreak p_{i})\}_{i=1}^{\infty}$\emph{ locally
converges} to a pointed metric space $\left(  \mathcal{X},d,x\right)  $ in the
sense of $C^{\infty}$--local submersions at $x$ if there is a Riemannian metric
$h$ on an open neighborhood $V\subset\mathbb{R}^{n}$ of $0$, a pseudogroup
$\Gamma$ of local isometries of $\left(  V,h\right)  $ such that the quotient
is well defined, an open set $U\subset\mathcal{X}$, and maps $\phi_{i}\co \left(
V,0\right)  \rightarrow\left(  \mathcal{M}_{i},p_{i}\right)  $ such that

\begin{enumerate}
\item $\left\{  \left(  \mathcal{M}_{i},d_{g_{i}},p_{i}\right)  \right\}
_{i=1}^{\infty}$ converges to $\left(  \mathcal{X},d,x\right)  $ in the
pointed Gromov--Hausdorff topology,

\item  the identity component of $\Gamma$ is a Lie group germ,

\item $\left(  V/\Gamma,\bar{d}_{h}\right)  $ is isometric to $\left(
U,d\right)  ,$

\item $\left(  \phi_{i}\right)  _{\ast}$ is nonsingular on $V$ for all
$i\in\mathbb{N},$ and

\item $h$ is the $C^{\infty}$ limit of $\phi_{i}^{\ast}g_{i}$ (uniform
convergence on compact sets together with all derivatives)
\end{enumerate}

where $\bar{d}_{h}$ is the induced distance in the quotient.
\end{definition}

We use the term convergence in $C^{\infty}$--local submersions since there is
actually a local submersion (actually a local covering) structure on manifolds
with bounded curvature defined by the exponential map. The convergence is
really of these structures. A local isometry $f$ of $V$ is a map from $V$ to
possibly a larger space which is a Riemannian isometry, i.e.\ $f^{\ast}h=h,$
usually for us it will be an injective map $B\left(  p,r/2\right)  \rightarrow
B\left(  p,r\right)  $ for some $p\in V$ and some $r>0.$ The global definition is:

\begin{definition}
A sequence of marked $n$--dimensional Riemannian manifolds $\{(\mathcal{M}%
_{i}^{n},\allowbreak g_{i},\allowbreak p_{i})\}_{i=1}^{\infty}$\emph{
converges} to a pointed metric space $\left(  \mathcal{X},d,x\right)  $ in the
sense of $C^{\infty}$--local submersions if for every $y\in\mathcal{X}$ there
exist $q_{i}\in\mathcal{M}_{i}$ such that $\left\{  \left(  \mathcal{M}%
_{i},g_{i},q_{i}\right)  \right\}  _{i=1}^{\infty}$\emph{ }locally converges
to $\left(  \mathcal{X},d,y\right)  $ in the sense of $C^{\infty}$--local
submersions at $y.$
\end{definition}

For the explanation of notation and definition of some of the terms, we refer
the reader to Section \ref{backgroundsection}. We are now ready to state our
main theorem. We refer to the metric $h$ on the ball $V$ corresponding to the
point $y\in\mathcal{X}$ as $h_{y}.$

\begin{theorem}
[Main Theorem]\label{C-infty local compactness}Let $C_{k}>0$ be constants for
$k\in\mathbb{N}$ and $\{(\mathcal{M}_{i}^{n},\allowbreak g_{i}\left(
t\right)  ,p_{i})\}_{i=1}^{\infty},$ where $t\in\left[  0,T\right]  ,$ be a
sequence of marked solutions to the Ricci flow on complete manifolds such
that
\[
\left|  Rm\left(  g_{i}\left(  t\right)  \right)  \right|  _{g_{i}\left(
t\right)  }\leq1
\]
and%
\[
\left|  D_{g_{i}\left(  t\right)  }^{k}Rm\left(  g_{i}\left(  t\right)
\right)  \right|  _{g_{i}\left(  t\right)  }\leq C_{k}%
\]
for all $i,k\in\mathbb{N}$ and $t\in\left[  0,T\right]  $. 

Then there is a
subsequence $\left\{  \left(  \mathcal{M}_{i_{k}},g_{i_{k}}\left(  t\right)
,p_{i_{k}}\right)  \right\}  _{k=1}^{\infty}$ and a one parameter family of
complete pointed metric spaces $\left(  \mathcal{X},d\left(  t\right)
,x\right)  $ such that for each $t\in\left[  0,T\right]  ,$ $\left(
\mathcal{M}_{i_{k}},d_{g_{i_{k}}\left(  t\right)  },p_{i_{k}}\right)  $
converges to $\left(  \mathcal{X}_{\infty},d_{\infty}\left(  t\right)
,x_{\infty}\right)  $ in the sense of $C^{\infty}$--local submersions and the
metrics $h_{y}\left(  t\right)  $ are solutions to the Ricci flow equation.
\end{theorem}

Peng Lu \cite{lu} has generalized Hamilton's compactness theorem in another
direction. He extended it to convergence of sequences of orbifold solutions to
the Ricci flow and showed that under similar conditions to Hamilton there is a
subsequence which converges to a orbifold solution to the Ricci flow. That
theorem uses an analogue of the injectivity radius bound. Although not treated
here, this work should be able to be extended to a larger class of sequences
such as orbifolds or potentially more general structures.

The primary advantage of the main theorem over other compactness theorems is
that it requires no control on the injectivity radius at all. Hence if we
consider blow ups about singularities as done in \cite{hamilton:singularities}%
, we do not have to worry about injectivity radius bounds, but may form
singularity models which collapse as well as ones which do not collapse. In
many cases it is believed that sequences arising from singularity blowups will
not collapse. Recent work of Perelman \cite{perelman} shows that sequences
arising from blow ups of finite time singularities on closed manifolds cannot
collapse. A study of collapsing singularity models may be useful in proving
other types of singularities, for instance Type IIb singularities, cannot
collapse. In \cite{chowglickensteinlu} Chow, Lu, and the author use the main
theorem to classify singularity models of three-dimensional manifolds arising
from sequences with almost nonnegative sectional curvature and diameters
tending to infinity. In addition, virtual limits are constructed to gain
geometric insight on the limits of sequences which collapse to non-manifold or
low dimensional limits. We hope to be able to use this theorem to better
understand the relationship between injectivity radius and the Ricci flow for
infinite time singularities.

We hope to be able to use the compactness theorem to understand long time
behavior of the Ricci flow analytically. Many nonsingular solutions collapse.
Isenberg--Jackson \cite{isenbergjackson} showed that some homogeneous metrics
collapse. Cara\-fora--Isenberg--Jackson \cite{carforaisenbergjackson},
Hamilton--Isenberg \cite{hamiltonisenberg}, Knopf \cite{knopf}, and
Knopf--McLeod \cite{knopfmcleod} demonstrated families of metrics which
quasiconverge to homogeneous solutions, which have the same asymptotic
behavior and hence collapse.

Perelman's work (\cite{perelman} and \cite{perelman2}) indicates a
decomposition of three-manifolds into ``thick'' and ``thin'' parts. The thin
parts are collapsing and can be understood topologically as a graph manifold.
We hope to be able to use a compactness theorem to understand the geometry of
these parts. A greater understanding of the manifolds or parts of the
manifolds may allow for the definition of a weak solution of the Ricci flow,
defined in terms of local submersions as in the compactness theorem. This
could be extremely useful in the study of manifolds of any dimension.

\rk{Acknowledgements}The author would like to thank his advisor Ben
Chow for suggesting this problem and for all of his help and
support. The author would also like to thank Nolan Wallach for all of
his help.

The author was partially supported by NSF grant DMS-0203926.

\section{Notation and Background\label{backgroundsection}}

For future use we would like to establish the following notation. A pointed
map written as $\phi\co \left(  X,x\right)  \rightarrow\left(  Y,y\right)  $
indicates that $\phi\left(  x\right)  =y.$ If $\left(  X,g_{X}\right)  $ and
$\left(  Y,g_{Y}\right)  $ are Riemannian manifolds, then $\phi$ is a local
isometry if $\phi^{\ast}g_{Y}=g_{X}.$ For any Riemannian manifold $\left(
M,g\right)  $ we use $d_{g}$ to denote the distance function on $M$ determined
by $g$ by integration along paths that makes $\left(  M,d_{g}\right)  $ into a
metric space. To avoid confusion, the term ``distance function'' shall be used
to refer to the metric for the metric space structure, and the term ``metric''
shall be reserved for Riemannian metrics. We also introduce the following
notation for balls in $\left(  X,d\right)  $%
\[
B_{d}\left(  p,r\right)  \doteqdot\left\{  q\in X:d\left(  p,q\right)
<r\right\}  .
\]
When it is clear which is the correct distance function, we may use
$B_{X}\left(  p,r\right)  $ instead. Given $Y\subset X$ we can define the ball
around $Y$ as
\[
B_{d}\left(  Y,r\right)  \doteqdot\left\{  q\in X:\exists y\in Y\text{ such
that }d\left(  y,q\right)  <r\right\}  .
\]
When no distance is specified, for instance $B\left(  0,r\right)  ,$ it is
assumed that we mean the Euclidean ball in $\mathbb{R}^{n}.$ We also shall
sometimes omit $0$ in the Euclidean ball and use $B\left(  r\right)  \doteqdot
B\left(  0,r\right)  .$

For the following we refer the reader to the treatment of local groups in
\cite{pontryagin}. We prefer the term pseudogroup since the term local group
usually refers to equivalence classes of arbitrarily small neighborhoods of
the identity. We shall actually need some of the components of our
pseudogroup, so group germs are not sufficient.

\begin{definition}
A topological space $G$ is a \emph{pseudogroup} if for certain pairs of
elements $a,b\in G$ there is a product $ab\in G$ such that

\begin{enumerate}
\item  If $ab,$ $\left(  ab\right)  c,$ $a\left(  bc\right)  $ are all
defined, then $\left(  ab\right)  c=a\left(  bc\right)  $

\item  If $ab$ is defined then for every neighborhood $W$ of $ab$ there exist
neighborhoods $U$ of $a$ and $V$ of $b$ such that $xy$ is defined for all
$x\in U$ and $y\in V$ and $xy\in W.$

\item  There is an identity element $e$ such that if $a\in G$ then $ae$ is
defined and $ae=a.$

\item  If for a pair of elements $a,b\in G,$ $ab$ is defined and $ab=e,$ then
$a$ is a left inverse for $b$ and we say $a=b^{-1}.$ If $b$ has a left
inverse, then for every neighborhood $U$ of $b^{-1}$ there exists a
neighborhood $V$ of $b$ such that every element $y\in V$ possesses a left
inverse $y^{-1}\in U.$
\end{enumerate}
\end{definition}

It is easy to show that if a pseudogroup $G$ acts as $B\left(  p,r/2\right)
\rightarrow B\left(  p,r\right)  $ as local isometries such that the relation
$x\sim\gamma x$ for all $\gamma\in G$ is an equivalence relation for $x\in
B\left(  p,r/4\right)  $ then the quotient $B\left(  p,r/4\right)  /G$ is
well-defined. When a pseudogroup acts in this way, we can define the push
forward of the distance function, denoted $\bar{d},$ which is a distance
function on the quotient. It is defined as%
\[
\bar{d}\left(  \left[  x\right]  ,\left[  y\right]  \right)  =\inf\left\{
d\left(  \gamma x,\gamma^{\prime}y\right)  :\gamma,\gamma^{\prime}\in
G\right\}
\]
where $\left[  x\right]  $ is the equivalence class of $x.$

We can also define Lie group germ as a special case of a pseudogroup.

\begin{definition}
A pseudogroup $G$ is a \emph{Lie group germ }if a neighborhood of the identity
of $G$ can be given a differentiable structure such that group multiplication
and inversion are differentiable maps when defined.
\end{definition}

Essentially this definition says that the pseudogroup $G$ is a subspace of a
Lie group which may not be a subgroup, i.e.\ closed under group operations.

Recall the definition of pointed Gromov--Hausdorff distance. Our definition is
the one used in \cite{fukayaboundary}. We first need the notion of an
$\epsilon$--pointed Gromov--Hausdorff approximation. Note that since a major
application is the blow up of a singularity and in that case the diameters may
go to infinity, we shall always work in the pointed category.

\begin{definition}
Let $\left(  X,d_{X}\right)  $ and $\left(  Y,d_{Y}\right)  $ be metric spaces
with $x_{0}\in X$ and $y_{0}\in Y$. \ A map $f\co (X,x_{0})\rightarrow(Y,y_{0})$
is an $\varepsilon$--pointed Gromov--Hausdorff approximation if it satisfies the following:

\begin{itemize}
\item $B_{Y}(y_{0},1/\varepsilon-\varepsilon)\subseteq B_{Y}(f\left[
B_{X}\left(  x_{0},1/\varepsilon\right)  \right]  ,\varepsilon)$
\item  For all $x,x^{\prime}\in B_{X}(x_{0},1/\varepsilon)$ we have
\[
\left|  d_{X}(x,x^{\prime})-d_{Y}\left(  f(x),f(x^{\prime})\right)  \right|
<\epsilon
\]
\end{itemize}
\end{definition}

Note that $f$ is not required to be continuous. \ The first condition
essentially says $f$ maps $X$ to almost all of $Y$, in the pointed sense, and
the second condition says that $f$ is approximately an isometry. Let $\vec{X}$
and $\vec{Y}$ respectively denote the pointed metric spaces $(X,d_{X},x_{0})$
and $(Y,d_{Y},y_{0}).$ We can now define the pointed Gromov--Hausdorff distance
$d_{GH}$ between $\vec{X}$ and $\vec{Y}$ as:%
\[
d_{GH}\left(  \vec{X},\vec{Y}\right)  \doteqdot\inf\left\{  \epsilon\,\left|
\begin{tabular}
[c]{c}%
$\exists\,\epsilon$--pointed Gromov--$\text{Hausdorff\thinspace\ }$\\
$\text{approximations\thinspace\ }f\co (X,x_{0})\rightarrow(Y,y_{0})\text{ }$\\
$\text{and }g\co (Y,y_{0})\rightarrow(X,x_{0})$%
\end{tabular}
\right.  \right\}  .
\]
We shall sometimes omit the dependence on the distance function if it is clear
what we mean. Gromov--Hausdorff distance is especially interesting since it can
compare manifolds with different topologies.

Note that it is not a distance in the usual sense, since the triangle
inequality does not hold exactly (see Proposition \ref{triangle
inequality}).  We have the following useful propositions most of whose
proofs we omit here.  For further details see
\cite{glickensteinthesis}. The first is an approximate triangle
inequality, which allows us to treat $d_{GH}$ essentially as a
distance.

\begin{proposition}
\label{triangle inequality}Let $\mathcal{\vec{X}}_{i}\doteqdot\left(
\mathcal{X}_{i},d_{i},O_{i}\right)  $ be pointed metric spaces for $i=1,2,3.$
If $d_{GH}\left[  \mathcal{\vec{X}}_{1},\mathcal{\vec{X}}_{2}\right]  \leq1/2$
and $d_{GH}\left[  \mathcal{\vec{X}}_{2},\mathcal{\vec{X}}_{3}\right]
\leq1/2$ then%
\[
d_{GH}\left[  \mathcal{\vec{X}}_{1},\mathcal{\vec{X}}_{3}\right]  \leq2\left(
d_{GH}\left[  \mathcal{\vec{X}}_{1},\mathcal{\vec{X}}_{2}\right]
,d_{GH}\left[  \mathcal{\vec{X}}_{2},\mathcal{\vec{X}}_{3}\right]  \right)  .
\]
\end{proposition}

The following is an associativity of convergence.

\begin{proposition}
\label{GH converge/ GH converge}Let $\mathcal{\vec{X}}_{i}\doteqdot\left(
\mathcal{X}_{i},d_{i},O_{i}\right)  $ and $\mathcal{\vec{X}}^{\prime}%
_{i}\doteqdot\left(  \mathcal{X}_{i}^{\prime},d_{i}^{\prime},O_{i}^{\prime
}\right)  $ be pointed metric spaces for $i\in\mathbb{N\cup}\left\{
\infty\right\}  .$ If $d_{GH}\left[  \mathcal{\vec{X}}_{i},\mathcal{\vec{X}%
}_{i}^{\prime}\right]  \leq\varepsilon$ for each $i\in\mathbb{N}$ and if
$\mathcal{\vec{X}}_{i}$ converges to $\mathcal{\vec{X}}_{\infty}%
\doteqdot\left(  \mathcal{X}_{\infty},d_{\infty},O_{\infty}\right)  $ and
$\mathcal{\vec{X}}_{i}^{\prime}$ converges to $\mathcal{\vec{X}}_{\infty
}^{\prime}\doteqdot\left(  \mathcal{X}_{\infty}^{\prime},d_{\infty}^{\prime
},O_{\infty}^{\prime}\right)  $ in GH then $d_{GH}\left[  \mathcal{\vec{X}%
}_{\infty},\mathcal{\vec{X}}_{\infty}^{\prime}\right]  \leq4\varepsilon.$
\end{proposition}

The following is an interaction between Gromov--Hausdorff convergence and
Lipschitz convergence.

\begin{proposition}
\label{GH converge/Lispchitz converge}Let $\left(  \mathcal{X}_{i},d_{i}%
,x_{i}\right)  $ and $\left(  \mathcal{X}_{i}^{\prime},d_{i}^{\prime}%
,x_{i}^{\prime}\right)  $ be pointed complete length spaces (for details see
\cite{buragoburagoivanov}) for $i\in\mathbb{N\cup}\left\{  \infty\right\}  .$
If for each $i\in\mathbb{N}$ there is a homeomorphism $\mathcal{X}%
_{i}\rightarrow\mathcal{X}_{i}^{\prime}$ taking $x_{i}$ to $x_{i}^{\prime}$
and Lipschitz with Lipschitz constant $c<\infty$ (independent of $i)$ and if
$\left(  \mathcal{X}_{i},d_{i},x_{i}\right)  $ converges to $\left(
\mathcal{X}_{\infty},d_{\infty},x_{\infty}\right)  $ and $\left(
\mathcal{X}_{i}^{\prime},d_{i}^{\prime},x_{i}^{\prime}\right)  $ converges to
$\left(  \mathcal{X}_{\infty}^{\prime},d_{\infty}^{\prime},x_{\infty}^{\prime
}\right)  $ in the pointed Gromov--Hausdorff topology then there is a Lipschitz
homeomorphism $\left(  \mathcal{X}_{\infty},x_{\infty}\right)  \rightarrow
\left(  \mathcal{X}_{\infty}^{\prime},,x_{\infty}^{\prime}\right)  $ with
Lipschitz constant $c.$
\end{proposition}

The proof of Proposition \ref{GH converge/Lispchitz converge} relies heavily
on the following fact.

\begin{proposition}
If $\mathcal{X}$ is a length space, then $\mathcal{X}$ is isometric to the
direct limit%
\[
\underset{\longrightarrow}{\lim}B_{\mathcal{X}}\left(  p,1/i\right)
\]
where the structure is just the obvious inclusions.
\end{proposition}

As a corollary we get the following important property of length spaces.

\begin{proposition}
\label{balls converge}If $\left(  \mathcal{X}_{i},d_{i},x_{i}\right)  $ are
length spaces converging to $\left(  \mathcal{X}_{\infty},d_{\infty}%
,x_{\infty}\right)  $ in the pointed Gromov--Hausdorff topology, then for every
$\rho>0,$ $B_{d_{i}}\left(  x_{i},\rho\right)  $ converge to $B_{d_{\infty}%
}\left(  x_{\infty},\rho\right)  $ in the pointed Gromov--Hausdorff topology.
\end{proposition}

The following relates equivalence of Riemannian metrics with Gromov--Haus\-dorff distance.

\begin{proposition}
\label{metrics close implies GH close}Let $g$ and $h$ be Riemannian metrics on
a smooth manifold $\mathcal{M}$ with $O\in\mathcal{M}$. Suppose there exist
$\delta>0$ such that
\[
\left(  1+\delta\right)  ^{-1}g\leq h\leq\left(  1+\delta\right)  g,
\]
then
\[
d_{GH}\left[  \left(  \mathcal{M},d_{g},O\right)  ,\left(  \mathcal{M}%
,d_{h},O\right)  \right]  <2\delta^{1/4}\left(  1+\delta\right)  ^{1/2}.
\]
\end{proposition}

Finally the following indicates how we can change basepoints.

\begin{proposition}
\label{limit points well defined}If $\left(  \mathcal{X}_{i},d_{i}%
,x_{i}\right)  $ converges to $\left(  \mathcal{X}_{\infty},d_{\infty
},x_{\infty}\right)  $ in the pointed\break Gromov--Hausdorff distance, then for
every $y_{\infty}\in\mathcal{X}_{\infty}$ there exist $y_{i}\in\mathcal{X}%
_{i}$ such that $\left(  \mathcal{X}_{i},d_{i},y_{i}\right)  $ converges to
$\left(  \mathcal{X}_{\infty},d_{\infty},y_{\infty}\right)  $ in the pointed
Gromov--Hausdorff distance.
\end{proposition}

Gromov proved the following precompactness result about Riemannian manifolds
with the Gromov--Hausdorff topology.

\begin{theorem}
{\rm\cite[Theorem 5.3]{gromovenglish}}\qua The set of $n$--dimensional pointed
Riemannian manifolds with Ricci curvature uniformly bounded below is
precompact with respect to the pointed Gromov--Hausdorff topology.
\end{theorem}

This says that every sequence has a subsequence which converges to a complete
metric space. The space is actually a length space with curvature bounded
below in the sense of Alexandrov (see \cite{buragoburagoivanov}).

In \cite{fukayaboundary} Fukaya made the definition of a smooth element of the
closure of the set of Riemannian manifolds, but since smooth has a different
meaning in our case, we instead use the term \emph{nice.} His definition is
the following.

\begin{definition}
A pointed metric space $(X,p)$ is called \emph{nice} if there exists

\begin{itemize}
\item  a neighborhood $U$ of $p$ in $X$

\item  a compact Lie group $G_{p}$, whose identity component is isomorphic to
a torus

\item  a faithful representation of $G_{p}$ into $O(m)$
\end{itemize}

such that $U$ is homeomorphic to $V/G_{p}$ for some neighborhood $V$ of 0 in
$\mathbb{R}^{m}$ (note that this is a linear action of a subgroup of the
orthogonal group on $\mathbb{R}^{m}$) and furthermore there is a $G_{p}%
$--invariant Riemannian metric $g$ such that $(V/G_{p},\bar{0})$ with distance
function $\bar{d}_{g}$ is isometric (as pointed metric spaces) to $(U,p)$.
(Here we used $\bar{0}$ to denote the equivalence class of 0 and $\bar{d}_{g}$
to denote the induced metric on the quotient space.)
\end{definition}

We note that $m$ may not equal $n$. In Theorem 0.5 Fukaya proved that nice
metric spaces are dense in the closure of $n$--dimensional Riemannian manifolds
with uniformly bounded curvature. Upon close inspection of the proof, however,
he proves the following theorem, not explicitly stated, which is exactly the
unparameterized version of our local theory.

\begin{theorem}
If $\left\{  \mathcal{M}_{i},g_{i},O_{i}\right\}  _{i\in\mathbb{N}}$ is a
sequence of $n$--dimensional pointed Riemannian manifolds such that
\[
\left|  Rm\left(  g_{i}\right)  \right|  _{g_{i}}\leq1
\]
and there exist $C_{k}>0$ such that
\[
\left|  D_{i}^{k}Rm\left(  g_{i}\right)  \right|  _{g_{i}}\leq C_{k}%
\]
then there exists a subsequence which converges to a metric space $\left(
\mathcal{X},d,x_{0}\right)  $ such that there exists a Riemannian metric $h$
on $B\left(  0,1\right)  \subset\mathbb{R}^{n},$ a pseudogroup $G_{\infty}$ of
local isometries acting as $B\left(  0,1/2\right)  \rightarrow B\left(
0,1\right)  $ such that $\left(  B_{\mathcal{X}}\left(  x_{0},1/4\right)
,d\right)  $ is isometric to $\left(  B\left(  0,1/4\right)  /G_{\infty}%
,\bar{d}_{h}\right)  .$
\end{theorem}

Furthermore there is some characterization of $G_{\infty}$; it is a local Lie
group with nilpotent Lie algebra and $G_{x_{0}}$ as in the definition of
``nice'' is the stabilizer of $0,$ i.e. the subgroup (not only
sub\emph{pseudo}group) $\left\{  \gamma\in G_{\infty}:\gamma\left(  0\right)
=0\right\}$.

We also need the following estimates of Hamilton.

\begin{lemma}
{\rm\cite[Lemma 2.4]{hamilton:compactness}}\label{hamilton bounds}\qua
Let $M$ be a Riemannian manifold with metric $G,K$ a compact subset of
$M$, and $G_{k}$ a collection of solutions to the Ricci flow defined
on neighborhoods of $K\times\lbrack\beta,\psi]$ with the time interval
$[\beta,\psi]$ containing $t=0$. Let $D$ denote covariant derivative
with respect to $G$ and $|\quad|$ the length of a tensor with respect
to $G$, while $D_{k}$ and $|\quad|_{k}$ are the same for
$G_{k}$. Suppose that

\begin{enumerate}
\item  the metrics $G_{k}$ are all uniformly equivalent to $G$ at $t=0$ on
$K$, so that
\[
cG(V,V)\leq G_{k}(V,V)\leq CG(V,V)
\]
for some constants $c>0$ and $C<\infty$ independent of $k$; and

\item  the covariant derivatives of the metrics $G_{k}$ with respect to the
metric $G$ are all uniformly bounded at $t=0$ on $K$, so that
\[
|D^{p}G_{k}|\leq C_{p}%
\]
for some constants $C_{p}<\infty$ independent of $k$ for all $p\geq1$; and

\item  the covariant derivatives of the curvature tensors $Rm_{k}$ of the
metrics $G_{k}$ are uniformly bounded with respect to the $G_{k}$ on
$K\times\lbrack\beta,\psi]$, so that
\[
|D_{k}^{p}Rm_{k}|_{k}\leq C_{p}^{\prime}%
\]
for some constants $C_{p}^{\prime}$ independent of $k$ for all $p\geq0.$
\end{enumerate}

Then the metrics $G_{k}$\ are uniformly bounded with respect to $G$\ on
$K\times\lbrack\beta,\psi]$, so that
\[
\tilde{c}G(V,V)\leq G_{k}(V,V)\leq\widetilde{C}G(V,V)
\]
for some constant $\tilde{c}$\ and $\widetilde{C}$\ independent of $k$, and
the covariant derivatives of the metrics $G_{k}$\ with respect to the metric
$G$\ are uniformly bounded on $K\times\lbrack\beta,\psi]$, so that
\[
|D^{p}G_{k}|\leq\widetilde{C}_{p}%
\]
for some constants $\widetilde{C}_{p}$\ independent of $k$\ with $\tilde
{c},\widetilde{C}$\ and $\widetilde{C}_{p}$\ depending only on $c,C,C_{p}%
$\ and $C_{p}^{\prime}$\ and the dimension.
\end{lemma}

\section{\smallskip Ricci flow estimates}

In this section we prove some estimates for ball sizes under evolution by the
Ricci flow. These are improvements of Hamilton's lemma given above as Lemma
\ref{hamilton bounds}. The ideas are wholly within Hamilton's work, but we
reproduce them here. Assume
\[
\left|  Rm\left(  g_{i}\left(  t\right)  \right)  \right|  _{g_{i}\left(
t\right)  }\leq C_{0}%
\]
and let $Rc_{i}\left(  t\right)  \doteqdot Rc\left(  g_{i}\left(  t\right)
\right)  .$ Let $V\in T_{p}\mathcal{M}$. We estimate $g_{i}\left(
t_{0}\right)  -g_{i}\left(  t\right)  $:%
\begin{align*}
\left[  g_{i}\left(  t\right)  -g_{i}\left(  t_{0}\right)  \right]  \left[
V,V\right]   &  =\int_{t_{0}}^{t}\frac{\partial}{\partial s}g_{i}\left(
s\right)  \left[  V,V\right]  ds\\
&  =\int_{t_{0}}^{t}-2Rc_{i}\left(  s\right)  \left[  V,V\right]  ds.
\end{align*}
and
\[
-\left|  Rc_{i}\left(  t\right)  \right|  g_{i}\left(  t\right)  \leq
Rc_{i}\left(  t\right)  \leq\left|  Rc_{i}\left(  t\right)  \right|
g_{i}\left(  t\right)
\]
And by our bounds on curvature,
\[
\left|  Rc_{i}\left(  t\right)  \left[  V,V\right]  \right|  \leq C_{0}%
g_{i}\left(  t\right)  \left[  V,V\right]  .
\]
Then we have
\[
\left|  \frac{\partial}{\partial t}\log\left(  g_{i}\left(  t\right)  \left[
V,V\right]  \right)  \right|  =2\frac{\left|  Rc_{i}\left(  t\right)  \left[
V,V\right]  \right|  }{g_{i}\left(  t\right)  \left[  V,V\right]  }\leq2C_{0}%
\]
so we can estimate%
\begin{align*}
\left|  \log\left(  g_{i}\left(  t\right)  \left[  V,V\right]  \right)
-\log\left(  g_{i}\left(  t_{0}\right)  \left[  V,V\right]  \right)  \right|
&  =\left|  \int_{t_{0}}^{t}\frac{\partial}{\partial t}\log\left(
g_{i}\left(  s\right)  \left[  V,V\right]  \right)  ds\right| \\
&  \leq\int_{t_{0}}^{t}\left|  \frac{\partial}{\partial t}\log\left(
g_{i}\left(  s\right)  \left[  V,V\right]  \right)  \right|  ds\\
&  \leq2C_{0}\left(  t-t_{0}\right)  .
\end{align*}
We conclude that%
\[
-2C_{0}\left(  t-t_{0}\right)  \leq\log\left(  g_{i}\left(  t\right)  \left[
V,V\right]  \right)  -\log\left(  g_{i}\left(  t_{0}\right)  \left[
V,V\right]  \right)  \leq2C_{0}\left(  t-t_{0}\right)
\]
which we can rewrite as
\[
e^{-2C_{0}\left(  t-t_{0}\right)  }\leq\frac{g_{i}\left(  t\right)  \left[
V,V\right]  }{g_{i}\left(  t_{0}\right)  \left[  V,V\right]  }\leq
e^{2C_{0}\left(  t-t_{0}\right)  }.
\]
Since $e^{x}=1+\delta\left(  x\right)  $ where $\delta\left(  x\right)
\rightarrow0$ as $x\rightarrow0,$ we also have
\[
\left|  g_{i}\left(  t\right)  \left[  V,V\right]  -g_{i}\left(  t_{0}\right)
\left[  V,V\right]  \right|  \leq\delta\left(  t-t_{0}\right)
\]
We can use this to estimate distances as well by integrating over paths and
get the following lemma.

\begin{lemma}
\label{metric bounds for RF}Suppose $\left(  \mathcal{M},g\left(  t\right)
\right)  $ are solutions to the Ricci flow for $t\in\left[  0,T\right)  $ such
that the sectional curvatures are uniformly bounded by $C_{0}.$ Then we have
the following:

\begin{enumerate}
\item  For all vectors $V$ and $t,t_{0}\in\left[  \alpha,\omega\right)  $%
\[
e^{-2C_{0}\left|  t-t_{0}\right|  }\leq\frac{g\left(  t\right)  \left[
V,V\right]  }{g\left(  t_{0}\right)  \left[  V,V\right]  }\leq e^{2C_{0}%
\left|  t-t_{0}\right|  }.
\]

\item  For all $\delta>0$ there exists an $\eta>0$ such that if $\left|
t-t_{0}\right|  <\eta$ then
\[
-\delta g\left(  t_{0}\right)  \leq g\left(  t\right)  -g\left(  t_{0}\right)
\leq\delta g\left(  t_{0}\right)  .
\]

\item  For all $\delta>0$ there exists an $\eta>0$ such that if $\left|
t-t_{0}\right|  <\eta$ then
\[
\left|  d_{g\left(  t\right)  }\left(  q,q^{\prime}\right)  -d_{g\left(
t_{0}\right)  }\left(  q,q^{\prime}\right)  \right|  \leq\delta^{1/2}%
d_{g\left(  t_{0}\right)  }\left(  q,q^{\prime}\right)
\]
for all $q,q^{\prime}\in\mathcal{M}.$
\end{enumerate}
\end{lemma}

\begin{proof}
Use the above arguments and let $\delta=\exp\left(  2C_{0}\left|
t-t_{0}\right|  \right)  -1$ and $\eta=\log\left(  \delta+1\right)  /\left(
2C_{0}\right)  .$
\end{proof}

As a corollary we can compute how the size of balls change with the Ricci flow.

\begin{proposition}
\label{containment}\bigskip If $g\left(  t\right)  $ is a solution to the
Ricci flow such that the curvature is bounded by the constant $1$ then for all
$\rho>0$
\begin{align*}
B_{g\left(  t\right)  }\left(  0,r\left(  t\right)  \rho\right)   &  \subset
B_{g\left(  0\right)  }\left(  0,\rho\right) \\
B_{g\left(  0\right)  }\left(  0,r\left(  t\right)  \rho\right)   &  \subset
B_{g\left(  t\right)  }\left(  0,\rho\right)
\end{align*}
where%
\[
r\left(  t\right)  =\frac{1}{1+\left(  e^{2t}-1\right)  ^{1/2}}.
\]
\end{proposition}

\begin{proof}
\bigskip Our estimates in Lemma \ref{metric bounds for RF} show%
\[
\left|  d_{g\left(  0\right)  }\left(  q,q^{\prime}\right)  -d_{g\left(
t\right)  }\left(  q,q^{\prime}\right)  \right|  \leq\delta\left(  t\right)
^{1/2}d_{g\left(  t\right)  }\left(  q,q^{\prime}\right)
\]
where $\delta\left(  t\right)  =\left(  e^{2t}-1\right)  ,$ so
\begin{align*}
d_{g\left(  0\right)  }\left(  p_{i},q\right)   &  \leq\left|  d_{g\left(
0\right)  }\left(  p_{i},q\right)  -d_{g\left(  t\right)  }\left(
p_{i},q\right)  \right|  +d_{g\left(  t\right)  }\left(  p_{i},q\right) \\
&  \leq\left(  1+\delta\left(  t\right)  ^{1/2}\right)  d_{g\left(  t\right)
}\left(  p_{i},q\right)
\end{align*}
so $B_{g\left(  t\right)  }\left(  0,\left(  1+\delta\left(  t\right)
^{1/2}\right)  ^{-1}\rho\right)  \subset B_{g\left(  0\right)  }\left(
0,\rho\right)  .$ Recall that we define%
\[
r\left(  t\right)  \doteqdot\frac{1}{1+\left(  e^{2t}-1\right)  ^{1/2}},
\]
so this says $B_{g\left(  t\right)  }\left(  0,r\left(  t\right)  \rho\right)
\subset B_{g\left(  0\right)  }\left(  0,\rho\right)  .$ Conversely, we have
\begin{align*}
d_{g\left(  t\right)  }\left(  p_{i},q\right)   &  \leq\left|  d_{g\left(
t\right)  }\left(  p_{i},q\right)  -d_{g\left(  0\right)  }\left(
p_{i},q\right)  \right|  +d_{g\left(  0\right)  }\left(  p_{i},q\right) \\
&  \leq\left(  1+\delta\left(  t\right)  ^{1/2}\right)  d_{g\left(  0\right)
}\left(  p_{i},q\right)
\end{align*}
so $B_{g\left(  0\right)  }\left(  0,r\left(  t\right)  \rho\right)  \subset
B_{g\left(  t\right)  }\left(  0,\rho\right)  .$
\end{proof}

\section{Convergence of metric spaces}

In this section we prove that the metric spaces converge. We will prove a
slightly more general theorem, since only the continuity in the $t$ variable
of solutions of the Ricci flow is used.

\begin{theorem}
\label{GH family compactness}Let $\left\{  \left(  \mathcal{M}_{i}%
,g_{i}\left(  t\right)  ,p_{i}\right)  \right\}  _{i=1}^{\infty},$ where
$t\in\left[  0,T\right]  ,$ be a sequence of pointed Riemannian manifolds of
dimension $n$ which is continuous in the $t$ variable in the following way:
for each $\delta>0$ there exists $\eta>0$ such that if $t_{0},t_{1}\in\left[
0,T\right]  $ satisfies $\left|  t_{0}-t_{1}\right|  <\eta$ then%
\[
\left(  1+\delta\right)  ^{-1}g_{i}\left(  t_{0}\right)  \leq g_{i}\left(
t_{1}\right)  \leq\left(  1+\delta\right)  g_{i}\left(  t_{0}\right)
\]
for all $i>0,$ and such that
\[
Rc\left(  g_{i}\left(  t\right)  \right)  \geq cg_{i}\left(  t\right)
\]
where $c$ does not depend on $t$ or $i$. 

Then there is a subsequence $\left\{
\left(  \mathcal{M}_{i_{k}},g_{i_{k}}\left(  t\right)  ,p_{i_{k}}\right)
\right\}  _{k=1}^{\infty}$ and a 1-parameter family of complete pointed metric
spaces $\left(  \mathcal{X}_{\infty}\left(  t\right)  ,d_{\infty}\left(
t\right)  ,x_{\infty}\left(  t\right)  \right)  $ such that for each
$t\in\left[  0,T\right]  $ the subsequence converges to $\left(
\mathcal{X}_{\infty}\left(  t\right)  ,d_{\infty}\left(  t\right)  ,x_{\infty
}\left(  t\right)  \right)  $ in the pointed Gromov-Hausdorff topology.
\end{theorem}

\begin{proof}
\strut For simplicity, let's define the notation $\mathcal{\vec{M}}_{i}\left(
t\right)  \doteqdot\left(  \mathcal{M}_{i},g_{i}\left(  t\right)
,p_{i}\right)  $ and $\mathcal{\vec{X}}_{\infty}\left(  t\right)
\doteqdot\left(  \mathcal{X}_{\infty}\left(  t\right)  ,d_{\infty}\left(
t\right)  ,x_{\infty}\left(  t\right)  \right)  .$ We shall take a number of
subsequences, each of which we shall continue to index by $i.$

Since the Ricci curvatures are bounded below for each time $t\in\left[
0,T\right]  ,$ at each time we can find a subsequence which converges to a
metric space. Instead, let $\mathcal{D}$ be a countable dense subset of
$\left[  0,T\right]  .$ For each $t_{0}\in\mathcal{D}$ we can find a
subsequence of $\mathcal{\vec{M}}_{i}\left(  t_{0}\right)  $ converging to a
metric space $\mathcal{\vec{X}}_{\infty}\left(  t_{0}\right)  .$ We can
diagonalize to find a subsequence which converges for all $t_{0}\in
\mathcal{D}$ since $\mathcal{D}$ is countable. We now show that the
subsequence (which we index also by $i$) is convergent for all $t\in\left[
0,T\right]  .$ For each $t\in\left[  0,T\right]  $ there is a subsequence
which converges to $\mathcal{\vec{X}}_{\infty}\left(  t\right)  .$ We claim
that $\mathcal{\vec{M}}_{i}\left(  t\right)  $ converges to $\mathcal{\vec{X}%
}_{\infty}\left(  t\right)  $ not just for the subsequence, but for all
$i\rightarrow\infty.$

Fix $\varepsilon>0$ and $t\in\left[  0,T\right]  .$ By Proposition
\ref{metrics close implies GH close} there is a $\delta>0$ such that if
$\left(  1+\delta\right)  ^{-1}g_{i}\left(  t_{0}\right)  \leq g_{i}\left(
t\right)  \leq\left(  1+\delta\right)  g_{i}\left(  t_{0}\right)  $ then
$d_{GH}\left(  \mathcal{\vec{M}}_{i}\left(  t_{0}\right)  ,\mathcal{\vec{M}%
}_{i}\left(  t\right)  \right)  <\varepsilon.$ This implies by Proposition
\ref{GH converge/ GH converge} that%
\[
d_{GH}\left(  \mathcal{\vec{X}}_{\infty}\left(  t_{0}\right)  ,\mathcal{\vec
{X}}_{\infty}\left(  t\right)  \right)  \leq4\varepsilon
\]
by taking the limit of the subsequence which converges for $t.$ By our
assumption, there exists $\eta>0$ such that if $\left|  t-t_{0}\right|  <\eta$
then
\[
\left(  1+\delta\right)  ^{-1}g_{i}\left(  t_{0}\right)  \leq g_{i}\left(
t\right)  \leq\left(  1+\delta\right)  g_{i}\left(  t_{0}\right)  .
\]
Since $\mathcal{D}$ is a dense set, we can choose $t_{0}\in\mathcal{D}$ such
that $\left|  t-t_{0}\right|  \leq\eta.$ Thus there exists $t_{0}%
\in\mathcal{D}$ such that
\[
d_{GH}\left(  \mathcal{\vec{M}}_{i}\left(  t_{0}\right)  ,\mathcal{\vec{M}%
}_{i}\left(  t\right)  \right)  <\varepsilon.
\]
Finally, since $\mathcal{\vec{M}}_{i}\left(  t_{0}\right)  $ converges to
$\mathcal{\vec{X}}_{\infty}\left(  t_{0}\right)  $ we can find $I$ large
enough so that
\[
d_{GH}\left(  \mathcal{\vec{M}}_{i}\left(  t_{0}\right)  ,\mathcal{\vec{X}%
}_{\infty}\left(  t_{0}\right)  \right)  <\varepsilon
\]
for $i>I.$ Hence, using our triangle inequality (Proposition \ref{triangle
inequality}) for Gromov-Hausdorff distance, we have for any $\varepsilon>0$
there is an $I>0$ such that if $i>I$ then
\[
d_{GH}\left(  \mathcal{\vec{M}}_{i}\left(  t\right)  ,\mathcal{\vec{X}%
}_{\infty}\left(  t\right)  \right)  \leq2\left(  2\left(  \varepsilon
+\varepsilon\right)  +4\varepsilon\right)  =16\varepsilon.
\]
Hence $\mathcal{\vec{M}}_{i}\left(  t\right)  \rightarrow\mathcal{\vec{X}%
}_{\infty}\left(  t\right)  $ as $i\rightarrow\infty.$
\end{proof}

Note that we haven't explicitly used solutions of the Ricci flow, although
solutions of the Ricci flow with bounded curvature would satisfy the premises
of the theorem, as shown below. Also note that we haven't really used the
upper curvature bound yet. That will be done in the next section when we look
at local properties of the limit.

We can now strengthen this result. Since $\left(  \mathcal{M}_{i},g_{i}\left(
t\right)  \right)  $ and $\left(  \mathcal{M}_{i},g_{i}\left(  0\right)
\right)  $ are homeomorphic by Lipschitz homeomorphisms, we can show that
$\mathcal{X}_{\infty}\left(  t\right)  $ is homeomorphic to $\mathcal{X}%
_{\infty}\left(  0\right)  .$

\begin{proposition}
\label{limits homeomorphic}If $\left|  \frac{\partial}{\partial t}g_{i}\left(
t\right)  \left[  V,V\right]  \right|  \leq Cg_{i}\left(  t\right)  \left[
V,V\right]  $ for each vector $V$ and for all $t\in\left[  0,T\right]  $ then
$\left(  \mathcal{X}_{\infty}\left(  t\right)  ,x_{\infty}\left(  t\right)
\right)  $ is homeomorphic to $\left(  \mathcal{X}_{\infty}\left(  0\right)
,x_{\infty}\left(  0\right)  \right)  $ for all $t\in\left[  0,T\right]  .$
\end{proposition}

\begin{proof}
An easy argument shows that the condition on the derivative of the metric
implies that
\[
e^{-C^{\prime}t}g_{i}\left(  0\right)  \leq g_{i}\left(  t\right)  \leq
e^{C^{\prime}t}g_{i}\left(  0\right)
\]
where $C^{\prime}$ is a constant depending only on $C$ and dimension (see
\cite[Section 8]{hamilton:singularities} for details). Hence for all
$i\in\mathbb{N}$ and $t\in\left[  0,T\right]  ,$ we must have
\[
e^{-C^{\prime}T}g_{i}\left(  0\right)  \leq g_{i}\left(  t\right)  \leq
e^{C^{\prime}T}g_{i}\left(  0\right)  .
\]
This implies that $\left(  \mathcal{M}_{i},g_{i}\left(  t\right)  \right)  $
and $\left(  \mathcal{M}_{i},g_{i}\left(  0\right)  \right)  $ are uniformly
Lipschitz homeomorphic via the identity map, where for each $i$ the identity
map is a Lipschitz homeomorphism. By Proposition \ref{GH converge/Lispchitz
converge} we must have that $\left(  \mathcal{X}_{\infty}\left(  t\right)
,d_{\infty}\left(  t\right)  ,x_{\infty}\left(  t\right)  \right)  $ is
Lipschitz homeomorphic to $\left(  \mathcal{X}_{\infty}\left(  0\right)
,d_{\infty}\left(  0\right)  ,x_{\infty}\left(  0\right)  \right)  ,$ and the
homeomorphism takes $x_{\infty}\left(  t\right)  $ to $x_{\infty}\left(
0\right)  .$
\end{proof}

\begin{corollary}
If $\left\{  \left(  \mathcal{M}_{i},g_{i}\left(  t\right)  ,O_{i}\right)
\right\}  _{i\in\mathbb{N}}$ is a sequence of complete solutions of the Ricci
flow for $t\in\left[  0,T\right]  $ such that the curvatures are uniformly
bounded then there is a subsequence which converges to a family of complete
metric spaces $\left(  \mathcal{X}_{\infty},d_{\infty}\left(  t\right)
,x_{\infty}\right)  $ in the pointed Gromov-Hausdorff topology.
\end{corollary}

\begin{proof}
The Ricci flow with uniformly bounded curvature satisfies the assumption of
Theorem \ref{GH family compactness} and Proposition \ref{limits homeomorphic}.
\end{proof}

\section{Local theory}

In this section we complete the proof of the main theorem. We need to
understand the local theory of the limit metric spaces. We proceed according
to Fukaya in \cite{fukayaboundary}.

Note that we need only look at a neighborhood of the basepoint because of
Proposition \ref{limit points well defined}. Fix frames $\mathcal{F}_{i}$ for
$T_{p_{i}}\mathcal{M}_{i}.$ Take $\phi_{i}=\exp_{p_{i}}\circ\mathcal{F}_{i},$
where the exponential map is the exponential map for the metric $g_{i}\left(
0\right)  ,$ independent of $t.$ Then $\phi_{i}\co B\left(  0,1\right)
\rightarrow\mathcal{M}_{i},$ where the ball is taken with respect to the
Euclidean distance. Since the $\phi_{i}$ are constant with respect to time and
the metrics $g_{i}$ satisfy the Ricci flow equation, it follows that the
metrics $\tilde{g}_{i}\left(  t\right)  \doteqdot\phi_{i}^{\ast}g_{i}\left(
t\right)  $ satisfy the Ricci flow equation. We use Proposition
\ref{containment}. Since $r\left(  t\right)  $ is a decreasing function,
$r\left(  T\right)  \leq r\left(  t\right)  $ and hence $B\left(  0,r\left(
T\right)  ^{2}/4\right)  \subset B_{\tilde{g}_{i}\left(  t\right)  }\left(
0,r\left(  T\right)  /4\right)  \subset B\left(  0,1/4\right)  $.

By Hamilton's lemma (our Lemma \ref{hamilton bounds}) all space derivatives of
curvature are bounded, and hence all space and time derivatives of the metric
are bounded by a fixed metric (say the Euclidean metric), so there is a
subsequence which converges in $C^{\infty}\left(  B\left(  0,1\right)
\times\left[  0,T\right]  \right)  $ to a family of Riemannian metrics
$\tilde{g}_{\infty}\left(  t\right)  .$ The time derivatives follow from the
Ricci flow equation and bounds on the curvature and its derivatives. Because
the convergence is $C^{\infty},$ $\tilde{g}_{\infty}\left(  t\right)  $
satisfies the Ricci flow equation.

We define the pseudogroups $\Gamma_{i}$ in the same way that Fukaya does, as
continuous maps that preserve $\phi_{i}$, i.e.%
\[
\Gamma_{i}=\left\{  \gamma\in C\left(  B\left(  0,\frac{1}{2}\right)
,B\left(  0,1\right)  \right)  :\phi_{i}\circ\gamma=\phi_{i}\right\}  .
\]
It is easy to check that $\Gamma_{i}$ are pseudogroups which act as
equivalence relations (so the quotient is well-defined) on $B\left(
0,1/4\right)  $. Also $\Gamma_{i}$ are independent of $t$ and are isometries
of $\left(  B\left(  0,r\left(  T\right)  ^{2}\right)  ,\tilde{g}_{i}\left(
t\right)  \right)  $ for each $t\in\left[  0,T\right]  $ (simply because the
metrics $\tilde{g}_{i}\left(  t\right)  $ are defined as pullbacks by
$\phi_{i},$ which does not depend on $t$). This may appear very restrictive
and possibly surprising, but it is because the metrics are pulled back by a
common local covering map, so $\tilde{g}_{i}$ are constructed to have these
symmetries. Furthermore, they converge to a limit pseudogroup $\Gamma_{\infty
}$ (see Fukaya \cite[p. 9]{fukayaboundary}).

At this point Fukaya shows that $\left(  B\left(  0,1/2\right)  ,\Gamma
_{i}\right)  $ converges to $\left(  B\left(  0,1/2\right)  ,\Gamma_{\infty
}\right)  $ in the equivariant Gromov--Hausdorff topology, and hence $B\left(
0,1/2\right)  /\Gamma_{i}$ converges $B\left(  0,1/2\right)  /\Gamma_{\infty}$
in the Gromov--Hausdorff topology by \cite{fukayaorbifolds}. 

Since $\left(
B\left(  O_{i},1/2\right)  \right)  $ converge in the Gromov Hausdorff
topology to $B\left(  O_{\infty},1/2\right)  $ and also $B\left(
0,1/2\right)  /\Gamma_{i}$ are isometric to $B\left(  O_{i},1/2\right)  ,$ by
Proposition \ref{GH converge/Lispchitz converge}\break $B\left(  0,1/2\right)
/\Gamma_{\infty}$ is isometric to $B\left(  O_{\infty},1/2\right)  .$ This is
precisely what happens when $t=0,$ however we need to work a little harder
when $t\neq0.$ The reason is twofold. Firstly the Euclidean ball $B\left(
0,r\right)  $ is not a metric ball for $\tilde{g}_{i}\left(  t\right)  $
unless $t=0$ since the $\phi_{i}$ is the exponential map for $g_{i}\left(
t\right)  .$ Secondly we must be careful, because $(B\left(  0,1/2\right)
,\tilde{g}_{i}\left(  t\right)  )/\Gamma_{i}$ may not be isometric to
$(B_{g_{i}\left(  t\right)  }\left(  O_{i},1/2\right)  ,g_{i}\left(  t\right)
)$ but to some other set. However, the Gromov--Hausdorff convergence only
really makes good sense on metric balls. Thus we need to be very careful with
which sets we use.

Now consider%
\[
\mathcal{F}_{i}\left(  t\right)  \doteqdot\left(  \phi_{i}\right)
^{-1}\left[  \overline{B_{g_{i}\left(  t\right)  }\left(  p_{i},\frac{r\left(
T\right)  }{4}\right)  }\right]  \subset\overline{B\left(  0,1/4\right)  }.
\]
The $\mathcal{F}_{i}\left(  t\right)  $ form a set of closed subsets of the
compact set $\overline{B\left(  0,1/4\right)  },$ and so for each $t$ there is
a subsequence which converges to $\mathcal{F}_{\infty}\left(  t\right)  ,$ a
closed subset of $\overline{B\left(  0,1/2\right)  }.$ This is because the
closed subsets of a compact set are compact in the Hausdorff topology. We
claim that we can find a subsequence such that $\mathcal{F}_{i}\left(
t\right)  \rightarrow\mathcal{F}_{\infty}\left(  t\right)  $ for all
$t\in\left[  0,T\right]  $ where $\mathcal{F}_{\infty}\left(  t\right)  $ are
some closed sets. We can find a subsequence which converges for all time by
the following proposition.

\begin{proposition}
There exist closed sets $\mathcal{F}_{\infty}\left(  t\right)  $ and a
subsequence such that for each $t\in\left[  0,T\right]  ,$ $\mathcal{F}%
_{i}\left(  t\right)  $ converges to $\mathcal{F}_{\infty}\left(  t\right)  $
in the Hausdorff distance.
\end{proposition}

\begin{proof}
Take a countable dense subset $\mathcal{D}\subset\left[  0,T\right]  .$ Since
for each $t$ there is a subsequence which converges, we can diagonalize so
that there is a subsequence which converges for all $t\in\mathcal{D}.$ We then
show that this sequence converges for all $t\in\left[  0,T\right]  .$ Clearly
we can define $\mathcal{F}_{\infty}\left(  t\right)  $ from the subsequence
which converges. Now since $B_{g_{i}\left(  t\right)  }\left(  p_{i}%
,\frac{r\left(  T\right)  }{4}\right)  $ is arbitrarily close to
$B_{g_{i}\left(  s\right)  }\left(  p_{i},\frac{r\left(  T\right)  }%
{4}\right)  $ for some $s\in\mathcal{D}$ (due to the facts that $\mathcal{D}$
is dense and we have estimates from Proposition \ref{metric bounds for RF}) we
get convergence for all $t.$
\end{proof}

An easy argument (see \cite{glickensteinthesis} for details) shows that
$\left(  \mathcal{F}_{i}\left(  t\right)  ,\Gamma_{i}\right)  $ converges to
$\left(  \mathcal{F}_{\infty}\left(  t\right)  ,\Gamma_{\infty}\right)  $ in
the equivariant Gromov--Hausdorff distance since both converge separately in
the Hausdorff distance. Thus $\mathcal{F}_{i}\left(  t\right)  /\Gamma_{i}$
converges to $\mathcal{F}_{\infty}\left(  t\right)  /\Gamma_{\infty}$ in the
Gromov--Hausdorff distance by \cite[Theorem 2.1]{fukayacollapse}. Since
$\overline{B_{g_{i}\left(  t\right)  }\left(  p_{i},\frac{r\left(  T\right)
}{4}\right)  }$ must converge to $\overline{B_{d_{\infty}\left(  t\right)
}\left(  x_{\infty},\frac{r\left(  T\right)  }{4}\right)  }$ by Proposition
\ref{balls converge} we can use Proposition \ref{GH converge/Lispchitz
converge} to show that $\mathcal{F}_{\infty}\left(  t\right)  /\Gamma_{\infty
}$ is isometric to $\overline{B_{d_{\infty}\left(  t\right)  }\left(
x_{\infty},\frac{r\left(  T\right)  }{4}\right)  }.$ Let's call that isometry
\[
\xi_{t}\co \mathcal{F}_{\infty}\left(  t\right)  /\Gamma_{\infty}\rightarrow
\overline{B_{d_{\infty}\left(  t\right)  }\left(  x_{\infty},\frac{r\left(
T\right)  }{4}\right)  }.
\]
We need only find the appropriate open set in $\mathcal{X}_{\infty}$ to
complete the proof of the main theorem.

Notice that $B\left(  0,\frac{r\left(  T\right)  ^{2}}{4}\right)
\subset\mathcal{F}_{i}\left(  t\right)  $ for all $i$ because
\[
\phi_{i}\left(  B\left(  0,\frac{r\left(  T\right)  ^{2}}{4}\right)  \right)
=B_{g_{i}\left(  0\right)  }\left(  p_{i},\frac{r\left(  T\right)  ^{2}}%
{4}\right)  \subset B_{g_{i}\left(  t\right)  }\left(  p_{i},\frac{r\left(
T\right)  }{4}\right)
\]
for all $i,$ so it follows that $B\left(  0,\frac{r\left(  T\right)  ^{2}}%
{4}\right)  \subset\mathcal{F}_{\infty}\left(  t\right)  .$ We claim that
\[
\xi_{t}\left[  \pi\left(  B\left(  0,\frac{r\left(  T\right)  ^{2}}{4}\right)
\right)  \right]  \subset B_{d_{\infty}\left(  t\right)  }\left(  x_{\infty
},\frac{r\left(  T\right)  }{4}\right)
\]
is open and hence is open in $\mathcal{X}_{\infty}$ where
\[
\pi\co B\left(  0,1/4\right)  \rightarrow B\left(  0,1/4\right)  /\Gamma_{\infty}%
\]
is the projection map. Let $y=\xi_{t}\left[  \pi\left(  \tilde{y}\right)
\right]  ,$ where $\tilde{y}\in B\left(  0,\frac{r\left(  T\right)  ^{2}}%
{4}\right)  .$ Since $\xi_{t}$ is an isometry and $\pi$ is a quotient by
isometries, we have that
\begin{align*}
d_{\infty}\left(  t\right)  \left[  x_{\infty},y\right]   &  =d_{\infty
}\left(  t\right)  \left[  \xi_{t}\left[  \pi\left(  0\right)  \right]
,\xi_{t}\left[  \pi\left(  \tilde{y}\right)  \right]  \right] \\
&  =\bar{d}_{\tilde{g}_{\infty}\left(  t\right)  }\left(  \pi\left(  0\right)
,\pi\left(  \tilde{y}\right)  \right) \\
&  \leq d_{g_{\infty}\left(  t\right)  }\left(  0,\tilde{y}\right) \\
&  <\frac{r\left(  T\right)  ^{2}}{4}\\
&  <\frac{r\left(  T\right)  }{4}%
\end{align*}
and hence $y\in B_{d\left(  t\right)  }\left(  x_{\infty},\frac{r\left(
T\right)  }{4}\right)  .$ Note, the distance on the quotient is defined as
\[
\bar{d}\left(  \pi\left(  x\right)  ,\pi\left(  y\right)  \right)
\doteqdot\inf_{\gamma,\gamma^{\prime}\in\Gamma_{\infty}}\left\{  d\left(
\gamma\cdot x,\gamma^{\prime}\cdot y\right)  \right\}
\]
which justifies the third step. Now, $\pi\left(  B\left(  0,\frac{r\left(
T\right)  ^{2}}{4}\right)  \right)  =B\left(  0,\frac{r\left(  T\right)  ^{2}%
}{4}\right)  /\Gamma_{\infty}$ is open in $\mathcal{F}_{\infty}\left(
t\right)  /\Gamma_{\infty}$ because%
\[
\pi^{-1}\left(  B\left(  0,\frac{r\left(  T\right)  ^{2}}{4}\right)
/\Gamma_{\infty}\right)  =\bigcup_{\gamma\in\Gamma_{\infty}}\gamma\left[
B\left(  0,\frac{r\left(  T\right)  ^{3}}{4}\right)  \right]
\]
which, since all $\gamma\in\Gamma_{\infty}$ are isometries, is the union of
open sets and hence open. Thus
\[
\mathcal{U}\left(  t\right)  \doteqdot\xi_{t}\left[  \pi\left(  B\left(
0,\frac{r\left(  T\right)  ^{2}}{4}\right)  \right)  \right]  \subset
B_{d_{\infty}\left(  t\right)  }\left(  x_{\infty},\frac{r\left(  T\right)
}{4}\right)
\]
is open in $B_{d_{\infty}\left(  t\right)  }\left(  x_{\infty},\frac{r\left(
T\right)  }{4}\right)  ,$ and hence open in $\mathcal{X}_{\infty}.$ Thus there
is an open set $\mathcal{U}\left(  t\right)  \subset\mathcal{X}_{\infty}$
which is isometric to $B\left(  0,\frac{r\left(  T\right)  ^{2}}{4}\right)
/\Gamma_{\infty}$ with the induced quotient metric of $\tilde{g}_{\infty
}\left(  t\right)  .$ This completes the proof of the main theorem.

As a corollary, we get Hamilton's compactness theorem. For simplicity, we
state it with curvature bound $1$ and injectivity radius bound $1$ at every
point. In fact, we could get Hamilton's full result using the
Cheeger--Gromov--Taylor \cite{cheegergromovtaylor} or Cheng--Li--Yau
\cite{chengliyau} result which bounds how fast the injectivity radius can fall
off with distance.

\begin{theorem}
Let $\left\{  \left(  \mathcal{M}_{i},g_{i}\left(  t\right)  ,p_{i}%
,\mathcal{F}_{i}\right)  \right\}  _{i=1}^{\infty}$ be a sequence of complete
marked Riemannian manifolds of dimension $n$ evolving by the Ricci flow such
that the curvatures are uniformly bounded by $1$ and the covariant derivatives
of the curvature tensor are uniformly bounded, and, in addition, the
injectivity radii are uniformly bounded below by $1.$ Then there exists a
subsequence which converges to a marked $C^{\infty}$ Riemannian manifold
$\left(  \mathcal{M}_{\infty},g_{\infty}\left(  t\right)  ,p_{\infty
},\mathcal{F}_{\infty}\right)  $ in $C^{\infty}$ on compact subsets.
\end{theorem}

\begin{proof}
[Proof (sketch)]We use the main theorem (Theorem \ref{C-infty local
compactness}), where the pseudogroups are trivial because the injectivity
radii at every point is bounded below by $1.$ We now have embeddings $B\left(
0,r\left(  T\right)  \right)  \rightarrow\mathcal{X}_{\infty}$ and we claim
that these form a $C^{\infty}$ structure such that the local Riemannian
metrics fit together to form a global metric.

For $y\in\mathcal{M}_{i}$ look at the map $\phi_{i,y}\co B\left(  0,r\left(
T\right)  \right)  \rightarrow B_{g_{i}}\left(  y,r\left(  T\right)  \right)
\subset\mathcal{M}_{i},$ induced by a frame and the exponential map; now
consider the overlaps maps $\phi_{i,y}^{-1}\circ\phi_{i,y^{\prime}}.$ These
are geodesic coordinates in a fixed metric, so by the curvature bounds we can
bound all derivatives of the maps uniformly in $i.$ Thus we can take a
subsequence so that these maps converge in $C^{\infty}$, where we take
countably many points $\left\{  y_{i,j}\right\}  \subset\mathcal{M}_{i}$ so
that
\[
\mathcal{M}_{i}=\bigcup_{j=1}^{\infty}\phi_{i,y_{i,j}}\left(  B\left(
0,r\left(  T\right)  \right)  \right)  .
\]
(We used Arzela--Ascoli to see that there exists a convergent subsequence.) The
limit maps are smooth transition functions, and thus give $\mathcal{X}$ a
smooth structure. Hence $\mathcal{X}$ is a manifold and we can call it
$\mathcal{M}_{\infty},$ and define $p_{\infty}\doteqdot x_{\infty}.$ By
continuity we see that the limit metrics form a tensor on all of
$\mathcal{M}_{\infty},$ and hence we have a global limit Riemannian metric
$g_{\infty}\left(  t\right)  $ for all $t$ and since it evolves by the Ricci
flow in every coordinate patch, $g_{\infty}\left(  t\right)  $ is a solution
to the Ricci flow.
\end{proof}

\section{Examples}

In this section we give two examples of collapsing sequences of Ricci flow.
The first gives an example where we do not converge to the Ricci flow on the
collapsed manifold. The second is a homogeneous example of what happens to a
particular nonsingular solution.

\begin{example}
$\mathcal{M}_{i}=\left(  S^{1}\times S^{1},dr^{2}+\frac{1}{i}f\left(
r\right)  ^{2}ds^{2}\right)  $
\end{example}

The Gauss curvature of $\mathcal{M}_{i}$ is
\[
K=-\frac{f^{\prime\prime}}{f}%
\]
and so the Ricci flow on $\mathcal{M}_{i}$ is
\[
\frac{\partial g_{i}}{\partial t}=2\frac{f^{\prime\prime}}{f}g_{i}.
\]
Notice that if we project to the first component, this is not stationary, but
clearly the initial metrics converge to the standard metric on $S^{1}.$ Were
this the Ricci flow, $S^{1}$ would be stationary, but the image of the Ricci
flow is changing.

In this case we can produce a global covering structure since the exponential
map is a covering map, and we see that the structure in the limit is simply
the group $\mathbb{R}$ acting on the second component of the universal cover
$\mathbb{R\times R}$ with metric $d\tilde{r}^{2}+f\left(  \tilde
{r}\operatorname{mod}2\pi\right)  ^{2}d\tilde{s}^{2}.$ Note that the
two-dimensional Ricci flow converges to a constant curvature metric by
\cite{hamilton:surfaces}, and hence the collapsed metric converges to a fixed
size metric, the image of $S^{1}\times S^{1}$ under the quotient map.

\begin{example}
Nonsingular convergence of Nil.
\end{example}

Consider the Ricci flow solution starting with a Nil metric on $\mathbb{R}%
^{3},$ say start with
\[
g_{0}=A\left(  dz-xdy\right)  ^{2}+Bdy^{2}+Cdx^{2}%
\]
and then evolve by the Ricci flow, and call this metric $g\left(  t\right)  $
at time $t$. Let this live on the torus, i.e.\ we quotient by the isometry
group generated by the isometries
\begin{align*}
\left(  x,y,z\right)   &  \rightarrow\left(  x+1,y,z\right) \\
\left(  x,y,z\right)   &  \rightarrow\left(  x,y+1,z\right) \\
\left(  x,y,z\right)   &  \rightarrow\left(  x,y,z+1\right)  .
\end{align*}
An easy calculation along the lines of \cite{isenbergjackson} gives the Ricci
flow solution for the associated ODE to be
\begin{align*}
A\left(  t\right)   &  =\left(  2t+C_{1}\right)  ^{-1/2}\\
B\left(  t\right)   &  =C_{2}\left(  2t+C_{1}\right)  ^{1/2}\\
C\left(  t\right)   &  =C_{3}\left(  2t+C_{1}\right)  ^{1/2}%
\end{align*}
for some constants $C_{k}.$ Note that the solution is nonsingular. Now we may
take the following sequence of solutions:%
\[
g_{i}\left(  t\right)  =g\left(  i+t\right)
\]
and consider $\mathcal{M}_{i}=\left(  S^{1}\times S^{1}\times S^{1}%
,g_{i}\left(  t\right)  \right)  $. It is easy to see that the sequence
converges to two dimensional Euclidean space $\mathbb{E}^{2},$ (noncompact)
and the limit of the covers have the metrics $dz^{2}+C_{2}dy^{2}+C_{3}dx^{2}.$
We show this by rescaling. Consider the map%
\[
\phi_{i}\left(  x,y,z\right)  =\left(  \frac{x}{\left(  2i+C_{1}\right)
^{1/4}},\frac{y}{\left(  2i+C_{1}\right)  ^{1/4}},\left(  2i+C_{1}\right)
^{1/4}z\right)
\]
which is a map%
\[
\mathbb{R}^{3}/\sim_{2}\,\rightarrow\mathbb{R}^{3}/\sim_{1}=S^{1}\times
S^{1}\times S^{1}%
\]
where equivalence $\sim_{1}$ comes from the isometries above and the
equivalence $\sim_{2}$ comes from the induced maps%
\begin{align*}
\left(  x,y,z\right)   &  \rightarrow\left(  x+\left(  2i+C_{1}\right)
^{1/2},y,z\right) \\
\left(  x,y,z\right)   &  \rightarrow\left(  x,y+\left(  2i+C_{1}\right)
^{1/2},z\right) \\
\left(  x,y,z\right)   &  \rightarrow\left(  x,y,z+\left(  2i+C_{1}\right)
^{-1/2}\right)  .
\end{align*}
These are the groups $\Gamma_{i}.$ It is clear that these group actions become
translation in the $z$--direction, which is the group $\Gamma_{\infty}$ in the
limit as $i\rightarrow\infty.$ Now consider
\begin{align*}
4\phi_{i}^{\ast}g_{i}\left(  t\right)   &  =\left(  2t+2i+C_{1}\right)
^{-1/2}\left[  \left(  2i+C_{1}\right)  ^{1/2}dz^{2}-2xdzdy+\left(
2i+C_{1}\right)  ^{-1}x^{2}dy^{2}\right] \\
&  \quad+C_{2}\frac{\left(  2t+2i+C_{1}\right)  ^{1/2}}{\left(  2i+C_{1}%
\right)  ^{1/2}}dy^{2}+C_{3}\frac{\left(  2t+2i+C_{1}\right)  ^{1/2}}{\left(
2i+C_{1}\right)  ^{1/2}}dx^{2}%
\end{align*}
so for $t=0,$ as $i\rightarrow\infty$ we get convergence to a Euclidean space
with metric $dz^{2}+C_{2}dy^{2}+C_{3}dx^{2}$ on the universal cover. Actually,
since
\[
\lim_{i\rightarrow\infty}\frac{\left(  2t+2i+C_{1}\right)  ^{1/2}}{\left(
2i+C_{1}\right)  ^{1/2}}=1
\]
we get that it converges to this for each time $t.,$ The manifolds converge in
Gromov--Hausdorff to the quotient, which is a steady state Euclidean plane
$\mathbb{E}^{2}$.

\end{document}